\documentclass{svjour3}

\usepackage{amsmath,amsfonts,amssymb}
\usepackage{graphicx}
\usepackage{bm}

\setkeys{Gin}{draft=false}
\usepackage{natbib}

\authorrunning{LEVEQUE, WAAGAN, GONZ\'ALEZ, RIM, and LIN}

\titlerunning{Generating Random Earthquake Events for PTHA}

\journalname{Submitted}

\newcommand{\ignore}[1]{}
\newcommand{\DB}{\Delta \bm B}
\def\reals{{{\rm l} \kern -.15em {\rm R} }}
\newcommand{\etamax}{\eta_{max}}
\newcommand{\DBshore}{\Delta B_{\scriptstyle{shore}}}
\newcommand{\dstrike}{d_{\scriptstyle{strike}}}
\newcommand{\ddip}{d_{\scriptstyle{dip}}}
\newcommand{\rstrike}{r_{\scriptstyle{strike}}}
\newcommand{\rdip}{r_{\scriptstyle{dip}}}
\newcommand{\ddepth}{d_{\scriptstyle{depth}}}
\newcommand{\bb}{{\bm b}}
\newcommand{\be}{{\bm e}}
\newcommand{\bs}{{\bm s}}
\newcommand{\bmu}{{\bm \mu}}
\newcommand{\bsigma}{\hat {\bm C}}
\newcommand{\bg}{{\bm g}}
\newcommand{\bv}{{\bm v}}
\newcommand{\bz}{{\bm z}}
\newcommand{\zm}{\bz^{[m]}}
\newcommand{\bzero}{{\bm 0}}
\newcommand{\bfC}{{\bm C}}
\newcommand{\bfI}{{\bm I}}
\newcommand{\bfV}{{\bm V}}
\newcommand{\bfLambda}{{\bm \Lambda}}
\newcommand{\bfTheta}{{\bm \Theta}}
\newcommand{\corr}{{\tt corr}}

\newcommand{\cref}[1]{(\ref{#1})}
\newcommand{\Cref}[1]{(\ref{#1})}
\newcommand{\Fig}[1]{Figure~\ref{#1}}
\newcommand{\Sec}[1]{Section~\ref{#1}}

\begin{document}

\title{Generating Random Earthquake Events for 
Probabilistic Tsunami Hazard Assessment}

\author{R. J. LeVeque, K. Waagan, F. I.
Gonz\'alez, D. Rim, and G. Lin}

\institute{%
Department of Applied Mathematics, University of Washington,
Seattle, WA. 
\and
Forsvarets Forskningsinstitutt, Oslo, Norway.
\and
Department of Earth and Space Sciences,
University of Washington, Seattle, WA.
\and 
Department of Applied Mathematics, University of Washington,
Seattle, WA.
\and
Department of Mathematics, Purdue University, West Lafayette, IN.
}

\maketitle

\begin{abstract}
In order to perform probabilistic tsunami hazard assessment (PTHA) based on
subduction zone earthquakes, it is
necessary to start with a catalog of possible future events along with the
annual probability of occurance, or a probability distribution of such events
that can be easily sampled.  For nearfield events, the distribution of
slip on the fault can have a significant effect on the resulting tsunami. We
present an approach to defining a probability distribution based on
subdividing the fault geometry into many subfaults and prescribing a desired
covariance matrix relating slip on one subfault to slip on any other
subfault.  The eigenvalues and eigenvectors of this matrix are then used to
define a Karhunen-Lo\`eve expansion for random slip patterns. This is similar
to a spectral representation of random slip based on Fourier series but
conforms to a general fault geometry.  We show that
only a few terms in this series are needed to represent the features of the
slip distribution that are most important in tsunami generation, first with a
simple one-dimensional example where slip varies only in the down-dip
direction and then on a portion of the Cascadia Subduction Zone.
\end{abstract}


\section{Introduction}
\label{sec:intro}

\begin{figure}
\hfil\includegraphics[width=.45\textwidth]{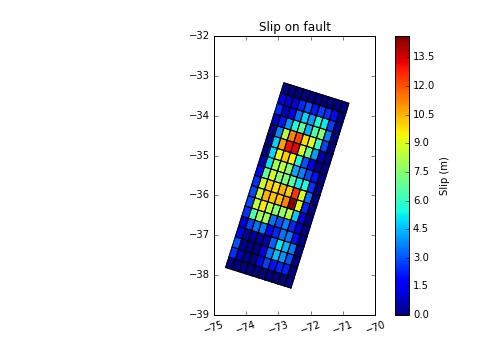}\hfil
\hfil\includegraphics[width=.45\textwidth]{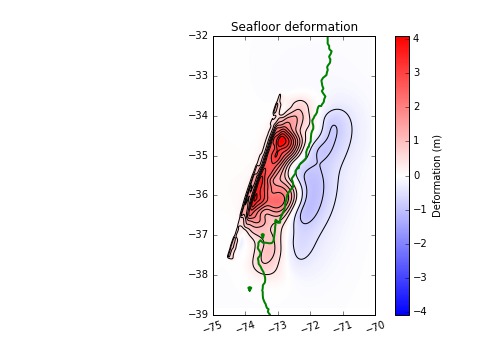}\hfil
\caption{\label{fig:chile} 
An example of slip distributed on a fault plane, from the USGS inversion 
of the 27 February 2010 event off Maule, Chile \cite{chile2010-usgs}. 
The plot on the right shows
the resulting sea floor deformation computed using the Okada model, with the
coast line in green.
  }
\end{figure}

Computer simulation of tsunamis resulting from subduction zone earthquakes can
be performed using a variety of available software packages, most of which
implement the two-dimensional shallow water equations and require the vertical 
seafloor motion resulting from the earthquake as the input to initiate the
waves.  For recent past events this can be approximated based on source
inversions; one example is shown in \Fig{fig:chile}.  
However, there are several situations in which it is desirable to
instead generate hypothetical future earthquakes. In particular, recent work on 
probabilistic tsunami hazard assessment (PTHA) has focused on producing maps
that indicate the annual probability of flooding exceeding various depths and
can provide much more information than a single ``worst considered case''
inundation map (for example \cite{
Adamstide,
GeistParsons2006,
GeistParsons:thesea,
Gonzalez:2009p452,
JaimesReinosoEtAl2016,
LovholtPedersenEtAl2012,
WitterZhangEtAl2013}).  
This requires running tsunami simulations for many potential earthquakes and
combining the results based on the annual probability of each, or using a Monte
Carlo approach to sample a presumed probability density of potential
earthquakes.  Generating a large number of hypothetical events can also be
useful for testing inversion methods that incorporate tsunami data, such as the
current DART buoy network, 
or that might give early tsunami warnings for the nearshore
(e.g., \cite{MelgarAllenEtAl2016}). 
Both Probabilistic Seismic Hazard Assessment (PSHA)
and PTHA are also fundamental tools in the
development of building codes that are critical in the design of structures
able to withstand seismic and tsunami forces (e.g., \cite{Chock2015}).

The primary goal of this paper is to introduce a general approach
to generating hypothetical earthquakes, by producing random slip
patterns on a pre-specified fault geometry.  Similar techniques have been used
in past studies, particularly for the generation of seismic waves in PSHA,
which has a longer history than PTHA.  
A variety of techniques have been proposed for generating random seismic
waveforms, see for example 
\cite{Anderson2015,DregerBerozaEtAl2015,Frankel1991,GhofraniAtkinsonEtAl2013,
GuatteriMaiEtAl2003,LavalleeLiuEtAl2006,MotazedianAtkinson2005,MaiBeroza2002}.
One approach is to use a spectral representation of the slip pattern as
a Fourier series with random coefficients that decay at a specified rate based
on the desired smoothness and correlation length of the slip patterns, e.g as
estimated from past events in the work of \cite{MaiBeroza2002}.
A random Fourier series can also be trimmed down to the desired
non-rectangular fault geometry, possibly with some tapering to zero
slip at some edges of the fault.  Different correlation lengths can
be specified in the strike and slip directions, if these directions
are used as the horizontal coordinates in the Fourier representation
and the fault is roughly rectangular.

\begin{figure}
\hfil\includegraphics[width=.9\textwidth]{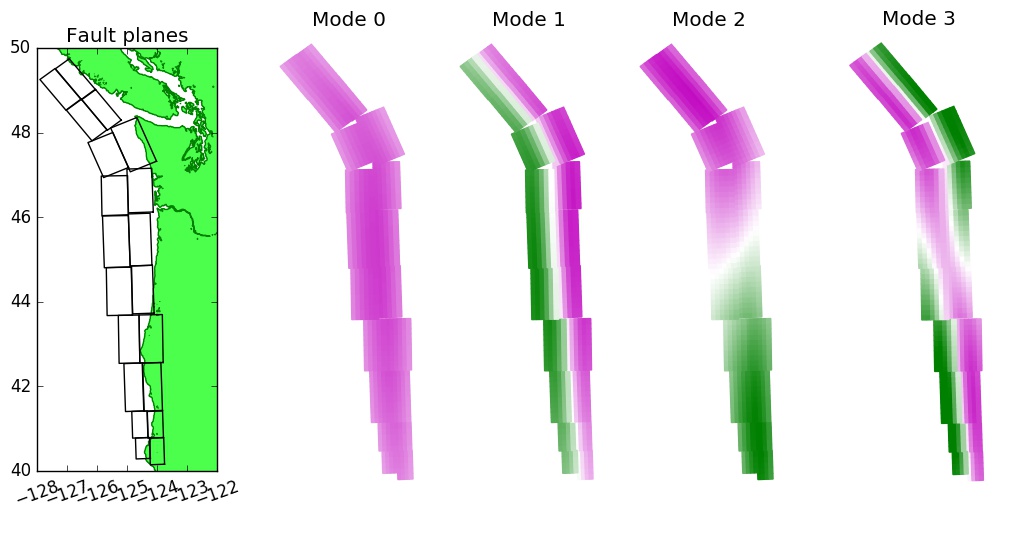}
\caption{\label{fig:csz1} 
Subdivision of the Cascadia Subduction Zone into 20 subfaults.  These are further
divided into 865 subfaults to compute the modes shown, which are the first
four eigenvectors of the $865\times 865$ covariance matrix as might be used in
a Karhunen-L\`oeve expansion.  Magenta and green are used to indicate positive
and negative entries in the eigenmodes.
  }
\end{figure}

Our approach is essentially the same on a rectangular fault
but generalizes easily to other fault geometries by using a Karhunen-Lo\`eve
expansion.  This work was
motivated in particular by the need to model events on the curving
Cascadia Subduction Zone (CSZ), which lies offshore North America
and runs nearly 1200 km from Northern California up to British
Columbia, see \Fig{fig:csz1}.

The fault is subdivided
into many rectangular subfaults and a value of the slip $s_i$ is
assigned to the $i$th subfault.  If here are $N$ subfaults, then
this defines a vector $\bs \in \reals^N$.  Initially we assume that the moment
magnitude $M_w$ of the earthquake (which depends on the total slip
summed over all subfaults) has been prescribed, and also that the
desired mean slip $\bmu \in \reals^N$ and covariance matrix $\bsigma \in
\reals^{N\times N}$ are known.
The mean slip is a vector with components $\mu_i = E[s_i]$, the
expected value of the slip on the $i$th subfault, and the $N\times
N$ covariance matrix $\bsigma$ has components $\hat C_{ij} =
E[(s_i-\mu_i)(s_j-\mu_j)]$, which can also be expressed as the outer product
$\bsigma = E[(\bs-\bmu)(\bs-\bmu)^T]$, where $T$ denotes transposing the
vector.

The {\em Karhunen-Lo\`eve (K-L) expansion} (e.g.
\cite{GhanemSpanos1991,HuangQuekEtAl2001,karhunen1947lineare,loeve:vol1,
schwab-todor-multipole})
is a standard approach to representing a random field as a linear
combination of eigenvectors of the presumed covariance matrix $\bsigma$.
If the matrix $\bsigma$ has
eigenvalues $\lambda_k$ (ordered with $\lambda_0 > \lambda_1 > \cdots > 0$)
and corresponding eigenvectors $\bv_k$, then the K-L expansion expresses
the slip vector $\bs$ as
\begin{equation}\label{KL1}
\bs = \bmu + \sum_{k=1}^N z_k \sqrt{\lambda_k} \bv_k.
\end{equation} 
where the $z_k$ are {\em independent} normally distributed random
numbers $z_k \sim {\cal N}(0,1)$ with mean 0 and standard deviation 1.
This is described in more detail in \Sec{sec:KL} where we explain why this
gives random slip patterns with the desired mean and covariance.  
This expansion makes it easy to generate an arbitrary number of realizations
using standard software to generate ${\cal N}(0,1)$ random numbers.  

\Fig{fig:csz1} shows an example of the
first four eigenmodes for the CSZ using this approach, where the $N$ components
of each eigenvector are represented on the fault geometry using a color map in
which magenta is positive and green is negative.  Note that Mode 0 is roughly
constant over the fault, so adding a multiple of this mode modifies the total
slip and hence the magnitude $M_w$.  On the other hand the other modes have
both positive and negative regions and so adding a multiple of any of these
tends to redistribute the slip (e.g., 
up-dip / down-dip with Mode 1 or between north and south with Mode 2).  
As with Fourier series, higher order eigenmodes are more oscillatory.

If the presumed correlation lengths are long and the covariance is a
sufficiently smooth function of the distance between subfaults,
then the eigenvalues $\lambda_k$ decay rapidly (there is little
high-frequency content) and so the K-L series can often be truncated
to only a few terms, greatly reducing the dimension of the stochastic
space that must be explored.

The K-L series approach could also be used to generate random slip patterns for
generating seismic waves, e.g. for performing PSHA or testing seismic inversion 
algorithms.  In this case high-frequency components of the slip are very
important and the K-L expansion may not decay so quickly.  However, for tsunami
modeling applications the slip pattern on the fault is only used to generate
the resulting seafloor deformation.  This is a smoothing operation that
suppresses high frequencies.  In this paper we also explore this effect and
show that truncating the expansion to only a few terms may be sufficient for
many tsunami applications.  Reducing the dimension of the stochastic space is
important for efficient application of many sampling techniques that could be
used for PTHA analysis.

In this paper we focus on explaining
the key ideas in the context of a one-dimensional fault model (with variation
in slip only in the down-dip direction) and a two-dimensional example using
the southern portion of the CSZ.  However, we do not claim to have used the
optimal parameters for modeling this particular fault.  
We also do not fully explore PTHA applications here, and for illustration
we use some quantities of interest related to a tsunami
that are easy to compute from a given slip realization, rather than
performing a full tsunami simulation for each.  This allows us to
explore the statistics obtained from a large number of realizations
(20,000) in order to illustrate some possible applications of this
approach and explore the effects of truncating the K-L series.
Work is underway to model the CSZ in a realistic manner and couple
this approach with a full tsunami model.

The K-L expansion as described above generates a Gaussian random field, in
which each subfault slip $s_i$ has a normal distribution with mean $\mu_i$ and
variance $\hat C_{ii}$ and together they have a joint normal distribution
with mean $\bmu$ and  covariance matrix $\bsigma$.  
A potential problem with this representation
is that when the variance is large it is possible for the slip $s_i$ to be
negative on some subfaults.  Since we assume the rake is constant (e.g. 90
degrees for a subduction thrust event), this would correspond to subfaults that
are slipping in the wrong direction.  The same issue arises with Fourier series
representations and can be dealt with by various means, for example by simply
setting the slip to zero any place it is negative (and then rescaling to
maintain the desired magnitude).  This naturally changes the statistics of the
resulting distributions.   

Another approach is to instead posit that the random
slip can be modeled by a joint lognormal distribution, for which the
probability of negative values is zero.  
Random slip patterns with a joint lognormal distribution can be generated by
using the K-L expansion to first compute a Gaussian field and then
exponentiating each component of the resulting vector to obtain the slip on
each subfault.  
By choosing the mean $\bmu^g$ and covariance matrix $\bsigma^g$
for the Gaussian field properly, the resulting lognormal will have the desired
mean $\bmu$ and $\bsigma$ for the slip.  This is discussed 
in \Sec{sec:lognormal} and used in the two-dimensional example in
\Sec{sec:def2d}.

\section{Expressing slip using a Karhunen-Lo\`eve expansion}\label{sec:KL}
If the earthquake fault is subdivided into $N$ small rectangular
subfaults, then a particular earthquake realization can be described by
specifying the slip on each subfault, i.e. by a vector $\bs \in \reals^N$
where $s_i$ is the slip on the $i$th subfault.  Note that we are assuming that only
the slip varies from one realization to another; the geometry 
and rake (direction of slip on each subfault) are fixed, and the slip is
instantaneous and not time-dependent.  These restrictions could be relaxed
at the expense of additional dimensions in our space of realizations.

Initially assume we wish to specify that the slip is a Gaussian
random field with desired mean slip $\bmu\in\reals^N$ and covariance
matrix $\bsigma\in\reals^{N\times N}$, which we write as $\bs \sim {\cal
N}(\bmu,\bsigma)$.  Then we compute the eigenvalues
$\lambda_k$ of $\bsigma$ and corresponding normalized eigenvectors $\bv_k$ so
that the matrix of eigenvectors $\bfV$ (with $k$th column $\bv_k$) and
diagonal matrix of eigenvalues $\bfLambda$ satisfy $\bsigma = \bfV\bfLambda
\bfV^T$. Note that the covariance matrix is symmetric postive definite,
so the eigenvalues are always positive real numbers and the
eigenvectors can be chosen to be orthonormal, $\bfV^{-1} = \bfV^T$.

Then the K-L expansion \cref{KL1} can be written in matrix-vector form as
\begin{equation}\label{KLV}
\bs = \bmu + \bfV\bfLambda^{1/2} \bz,
\end{equation} 
where $\bz\in\reals^N$ is a vector of independent identically distributed ${\cal
N}(0,1)$ random numbers.
Realizations generated via the K-L expansion have the right statistics 
since we can easily compute that
$E[\bs] = \bmu$ (since $E[\bz] = \bzero$) and 
\begin{equation}\label{KLsigma}
\begin{split} 
E[(\bs-\bmu)(\bs-\bmu)^T] &= E[\bfV\bfLambda^{1/2}\bz\bz^T \bfLambda^{1/2}\bfV^T] \\
&= \bfV\bfLambda^{1/2} E[\bz\bz^T]  \bfLambda^{1/2}\bfV^T \\
&= \bfV\bfLambda \bfV^T = \bsigma
\end{split} 
\end{equation} 
using the fact that $\bfV$ and $\bfLambda$ are fixed and $E[\bz \bz^T] = \bfI$.
Note that the $\bz$ could be chosen from a different probability density with 
mean $\bzero$ and covariance matrix $\bfI$ 
and achieve the same covariance matrix $\bsigma$ with the K-L expansion,
although the $\bs$ would not have a joint normal distribution in this case.

\section{One-dimensional case: down-dip variation}\label{sec:1d}

We first illustrate this technique on
a simplified case, a rectangular fault plane
that is essentially infinitely long in the strike direction and with uniform
slip in that direction, similar to the test case used by
\cite{LovholtPedersenEtAl2012}.  The slip will only vary in the
down-dip direction, reducing the problem to a single space dimension.  
The fault width is 100 km, a typical width for subduction zone faults, and
is assumed to dip at $13^\circ$ from horizontal, with the upper edge at a
depth of 5 km below the sea floor. 

\begin{figure}
\hfil\includegraphics[width=.45\textwidth]{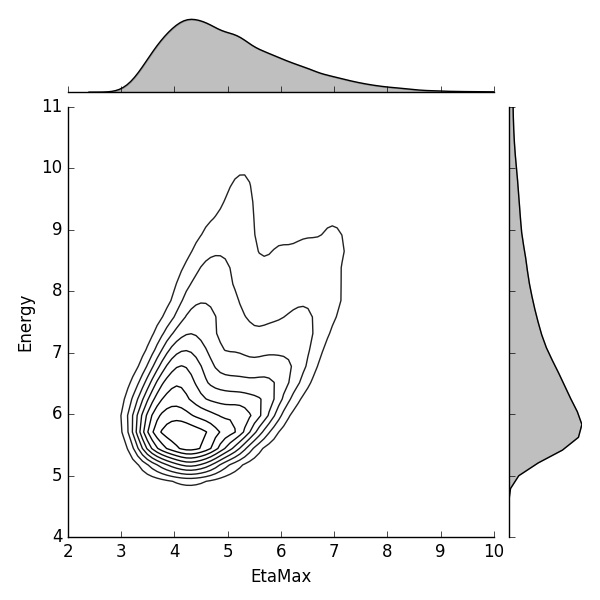}\hfil
\hfil\includegraphics[width=.45\textwidth]{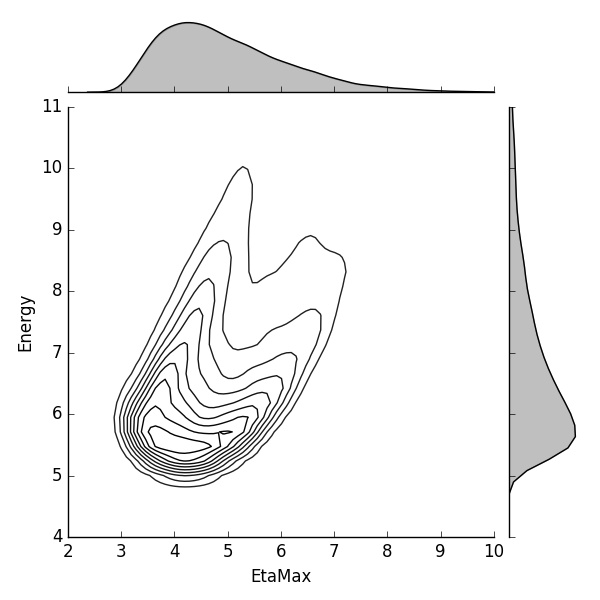}\hfil

\hfil\includegraphics[width=.45\textwidth]{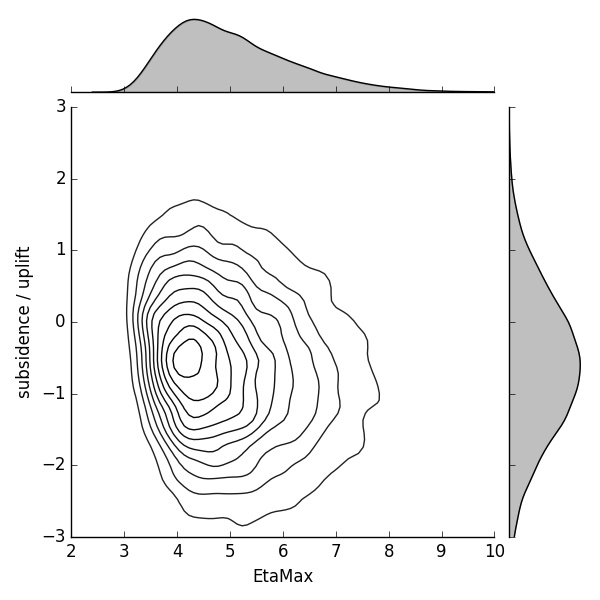}\hfil
\hfil\includegraphics[width=.45\textwidth]{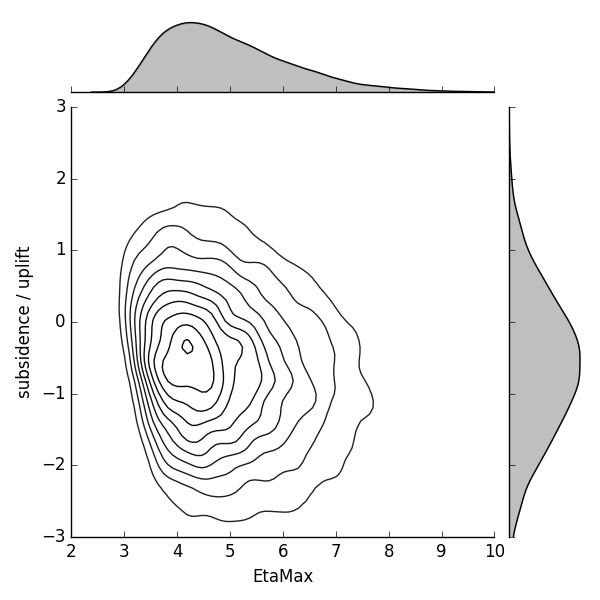}\hfil

\caption{\label{fig:joint} 
Joint and marginal probability densities for different quantities, comparing
the densities estimated using the 20-term expansion (left column) and the
3-term expansion (right column).  
The top row shows the joint density of $\etamax$ with
potential energy $E$ of the tsunami generated.
The bottom row shows the joint density
of $\etamax$ with $\DBshore$, the vertical displacement at the shore.
  }
\end{figure}

For the tests we perform here, we will focus on events of a single specified
magnitude. 
The moment magnitude $M_w$ is a function of the total slip integrated over the
entire fault plane, and also depends on the rigidity of the rock.  For
typical rigidity parameters, an average of 10 m of slip distributed over a
fault that is 100 km wide and 1000 km long would result 
in a magnitude $M_w \approx 9.0$
and so we fix the total slip to have this average.  If the fault were only
half as long, 500 km, then this would be a $M_w \approx 8.8$ event and
20 m average slip would be required for a
magnitude 9 event.  With the exception of the potential energy shown in
\Fig{fig:joint}, the quantities of interest considered
in this paper are all linear in the total
slip, however, so it does not really matter what value we choose.%
\footnote{
We follow \url{http://earthquake.usgs.gov/aboutus/docs/020204mag_policy.php}
and use $M_w = \frac 2 3 (\log_{10}(M_o) - 9.05)$ where the seismic moment
$M_o = $length $\times$ width$\times$(average slip)$\times$(rigidity) and 
set the rigidity to 
$3.55\times 10^{10}$ N-m for this calculation.}

\ignore{
Because of the
logarithmic dependence of $M_w$ on total slip, each doubling of
the slip would cause an increase of $M_w$ of about 0.2.
}

An important aspect of PTHA analysis is to consider possible events of
differing magnitudes as well, and take into account their relative
probabilities.  For smaller earthquakes, the
Gutenberg-Ricter relation approximately describes their relative frequency,
but for large subduction zone events that may have a recurrence time of
hundreds of years, there is generally no simple model for the variation of
annual probability with magnitude.  There may be a continuous distribution
of magnitudes or there may be certain ``characteristic earthquakes'' that
happen repeatedly after sufficient stress has built up.  The lack of
detailed data for past events over a long time period makes this difficult
to assess.

For the purposes of this paper, we assume that an earthquake of a particular
magnitude occurs and we wish to model the range of possible tsunamis that
can arise from such an event.  We thus discuss the relative probability of
different slip patterns and tsunamis given that an event of this magnitude
occurs, and so the probability density should integrate to 1.
This could then be used as one component in a full PTHA analysis by
weighting these results by the probability that an event of this magnitude
occurs and combining with similar results for other magnitudes.
Alternatively, one could introduce the magnitude as an additional stochastic
dimension and assume some probability density function for this, as we also
discuss further below.

We use $x$ to denote the
distance down-dip and split the fault into $N$ segments of equal
width $\Delta x$, where $N\Delta x$ is the total width of the fault in the dip
direction.  We then specify $N$ slips $s_i$ for $i=1,~2,~\ldots,~N$.  In our
one-dimensional experiments we take $N=200$. This is finer than one would
use in two dimensions and much finer than is needed to represent slip at an
adequate level for either seismic or tsunami modeling.  (For example, note from
\Fig{fig:chile} that the seismic inversion for this event represents the
slip as piecewise constant on a $18\times 10$ grid with only 10 segments in
the down-dip direction.)  One could certainly reduce the dimension of the
stochastic space below $N=200$ by using fewer subfaults.  However, we will
show that the dimension can be drastically reduced by instead using the K-L
expansion reduced to only a few terms (e.g. 3, for this one-dimensional
model the parameter choices below).  By starting with a fine discretization
of the fault, the eigenmodes used are smooth and perhaps better represent actual
slip patterns than piecewise constant functions over large subfaults.

\begin{figure}
\hfil\includegraphics[height=1.3in]{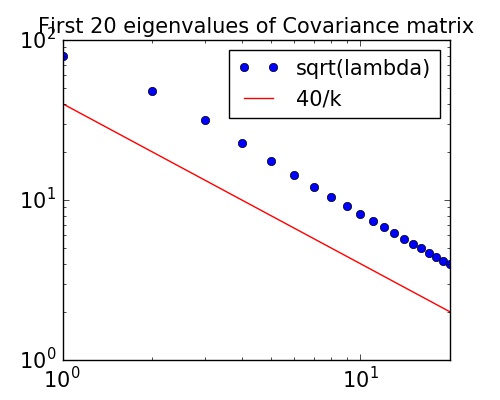}\hfil
\hfil\includegraphics[height=1.3in]{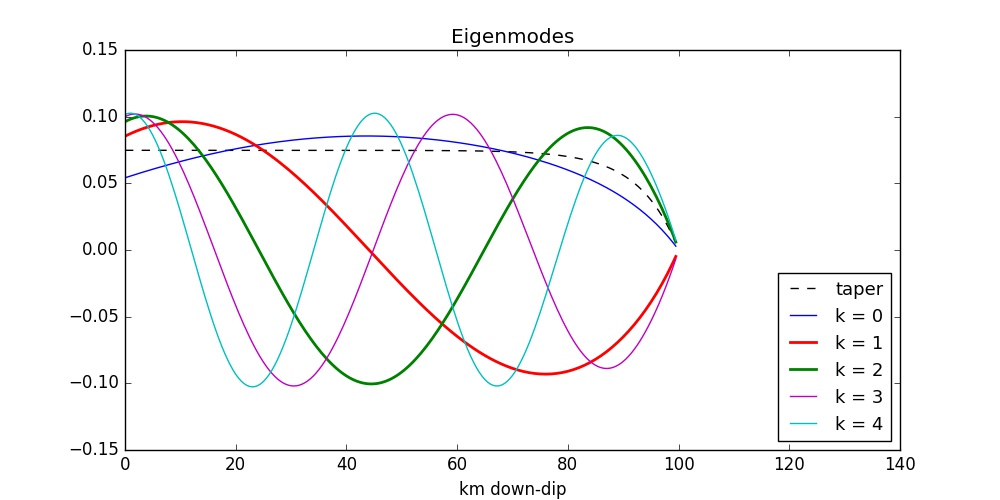}\hfil
\caption{\label{fig:modes1d} 
Eigenvalues decay like $1/k^2$ when the exponential autocorrelation function
is used.  The corresponding eigenvectors are similar to Fourier modes,
shaped by the taper.
  }
\end{figure}

We assume that the $N$ slips
are to be chosen randomly from a joint normal distribution
with mean $\bmu = [\mu_1,~\mu_2,~\ldots,~\mu_N]^T$.  The mean is chosen to be 
the desired taper, scaled to have the desired total slip.  
As an illustration of taper we use the function
\begin{equation}\label{taper}
\tau(d) = 1 - \exp(-20(d-d_{max})/d_{max})
\end{equation} 
where $d$ is the depth of a subfault and $d_{max}=22500$m  is the maximum depth of
the fault.  This function is close to 1 over most of the fault but tapers
toward the down-dip edge.
This taper, after scaling to give the mean slip,
is shown as the dashed line in \Fig{fig:modes1d}.  
Other tapers can be used instead, e.g. the taper propsed by \cite{WangHe2008}.

We set the desired covariance matrix to be 
$\hat C_{ij} = \sigma_i\sigma_j C_{ij}$ where $\sigma_i = \alpha\mu_i$ 
for some scalar
$\alpha\in\reals$ and $\bfC$ is the desired correlation matrix.
We take $\alpha = 0.75$, which tends to keep the slip positive everywhere, as
desired, while still giving reasonable variation in slip patterns.
The correlation matix is given by $C_{ij} =
\corr(|x_i-x_j|)$ in terms of some autocorrelation function (ACF) 
$\corr(r)$, and we choose
\begin{equation} \label{ACFexp}
\corr(r) = \exp(-r/r_0),
\end{equation} 
where the correlation length is set to $r_0 = 0.4W = 40$ km, i.e., 40\% of the
fault width as suggested by \cite{MaiBeroza2002}.

\Fig{fig:modes1d} shows the taper along with the first several
eigenvectors of the covariance matrix $\bfC$ (ordered based on the magnitude of
the eigenvalues, with eigenvectors normalized to have vector 2-norm equal to 1).  
Note that the lowest mode 0 looks very similar to the
taper.  Adding in a multiple of this mode will modify the total slip and
hence the magnitude, so we drop this mode from the sum.  The higher modes
are orthogonal to mode 0 and hence do not tend to change the total slip.
They look like Fourier modes that have been damped near the down-dip
boundary by the taper.

To create a random realization, we choose a vector $\bz$ of $N$ i.i.d.
Gaussian ${\cal N}(0,1)$ values $z_k$ for $k=0,~1,~\ldots,~N-1$.  
If we neglect the 0-mode and
truncate the expansion after $m$ terms, then this amounts to setting $z_0=0$
and $z_k=0$ for $k>m$.  We will denote such a $\bz$ vector by $\bz^{[m]}$.
The slip pattern can then be written as 
\begin{equation}\label{sm}
\bs = \bmu + \bfV \bfLambda^{1/2} \zm.
\end{equation} 

\begin{figure}[t]
\hfil\includegraphics[width=.9\textwidth]{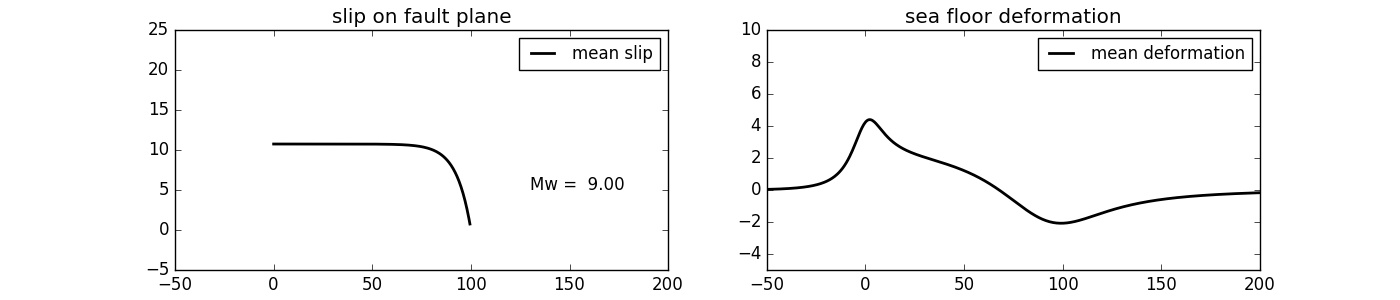}\hfil

\hfil\includegraphics[width=.9\textwidth]{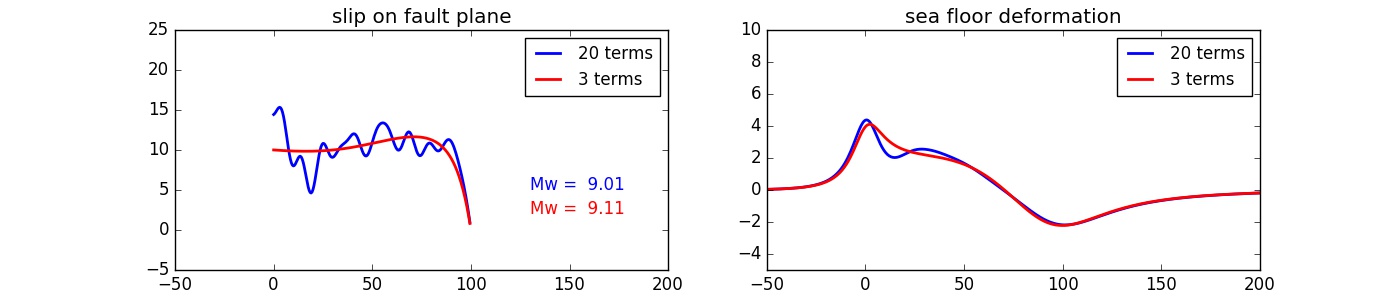}\hfil

\hfil\includegraphics[width=.9\textwidth]{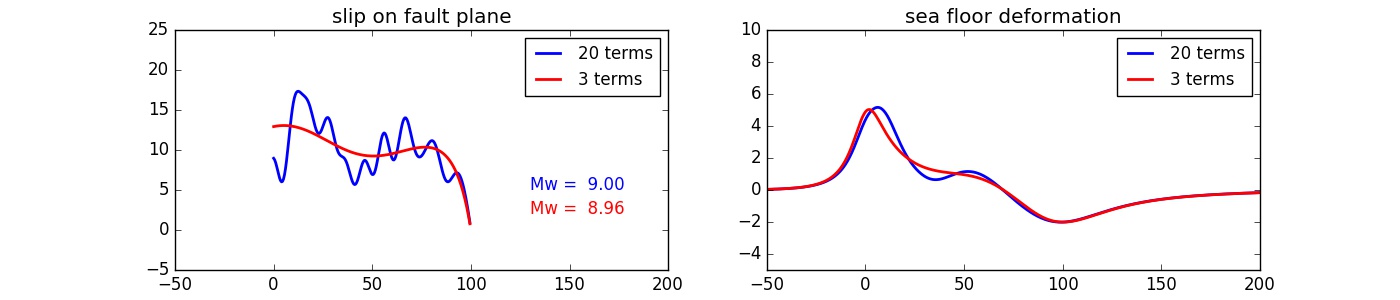}\hfil

\hfil\includegraphics[width=.9\textwidth]{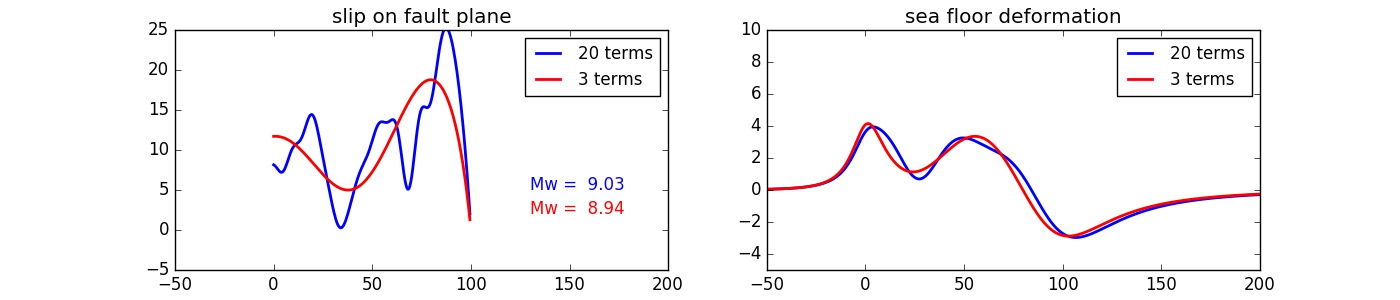}\hfil

\hfil\includegraphics[width=.9\textwidth]{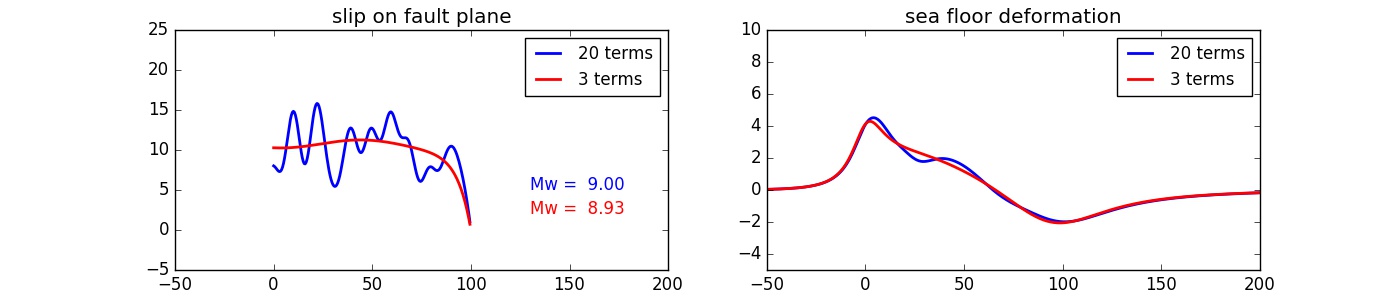}\hfil

\hfil\includegraphics[width=.9\textwidth]{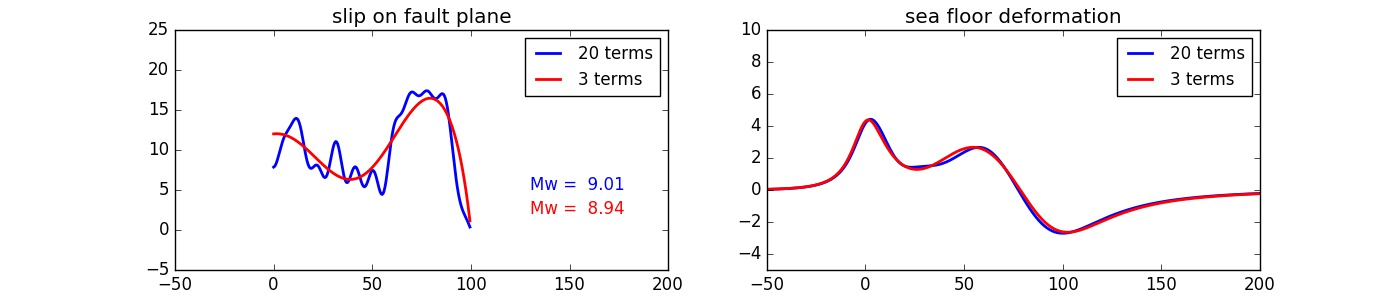}\hfil

\caption{\label{fig:realiz1} 
The left column shows slip on the fault plane of width $W = 100$ km.  The
right column shows the resulting seafloor deformation if the up-dip edge of
the fault plane is 5 km below the surface and it dips at $13^\circ$. The top
row shows the mean slip and resulting deformation.  The remaining rows show
random realizations using 20 terms  of a K-L expansion (blue) and the same
sum truncated to 3 terms (red). 
  }
\end{figure}

The left column of \Fig{fig:realiz1} shows the mean slip in the top plot, 
followed by 
several random realizations generated by the K-L expansion using 20 terms,
with the $\bz^{[20]}$ coefficients chosen as i.i.d. ${\cal N}(0,1)$ values.  
These are the blue curves in each plot.
In each
case, the slip is also shown when only 3 terms in the series are used (i.e.
$\bz^{[3]}$ is computed by leaving
$z_1,~z_2,~z_3$ unchanged from $\bz^{[20]}$ but with the 
higher terms dropped, equivalent to truncating the expansion at an earlier
point).  
These slip patterns, shown in red, are smoothed versions of
the 20-term slip patterns since the higher wave number components have been
surpressed.  In many cases there appears to be quite a large difference
between the 3-term and 20-term slips.  This is a reflection of the fact that
the eigenvalues do not decay all that quickly in this case.  There would be
faster decay if a longer correlation length were chosen, and much more
rapidly if the Gaussian autocorrelation function
were chosen instead of the exponential. 

In spite of the differences in the slip patterns, for tsunami modeling the
3-term series may still be adequate.
The right column of \Fig{fig:realiz1} shows the
sea floor deformations that result from the slips shown on the left.  
The blue curves show the deformation due to the 20-term sum while the red
curves show the deformation resulting from the truncated 3-term sum.  As
mentioned earlier, high wavenumber oscillations in the slip pattern are
filtered out by the Okada model used to compute
the seafloor deformation $\DB$ from the slip $\bs$.  If the $\DB$
are sufficiently similar
between the 3-term and the full K-L expansion, then there is no reason to
use more terms.  In this case we have reduced the stochastic space that
needs to be explored down to 3 dimensions.  There is much greater similarity for
some realizations than for others, and so below we examine the statistical
properties of this approximation by using a sample of 20,000 realizations.

Note that using only 3 modes may be too few 
for this particular set of fault parameters --- 
the comparisons shown in \Fig{fig:realiz1} would look more similar if a few
more terms were retained ---  but we will see that good statistical properties
are obtained even with this severe truncation.  How many terms are required
depends on various factors: not only the correlation structure of the
slip as discussed above, but also the depth of the fault plane.  The deeper
it is, the more smoothing takes place in the Okada model.  
Here we placed the top of the fault plane at 5 km depth. This is a typical
depth but there is variation between subduction zones.

To examine statistical properties of the 20-term sum and the 3-term
approximation, we generate 20,000 samples of each and compare some
quantities that are cheap to compute but that are important indicators of
the severity of the resulting tsunami.  Running a full tsunami model based
on the shallow water equations is not feasible for this large number of
realizations, but the quantities we consider will stand in as proxies for
the quantities one might actually want to compute, such as the maximum depth
of flooding at particular points onshore.  Moreover the distribution of
these proxy values can be used in a later stage to help choose particular
earthquake realizations for which the full tsunami model will be run. It is
desirable to run the model with judiciously chosen realizations for which
the proxy values are well distributed over the range of possible values.
The computed densities of the proxy values can also be used to weight the results
from the full model runs to accurately reflect the probabities of such
events.  This will be explored in detail in a future paper.

Computations for a large number of realizations can be sped up substantially
by realizing that the Okada model is linear in the slip, i.e. if the slip
vector is given by $\bs$ then the resulting sea floor deformation can be
written as $\DB = \bfTheta \bs$ for a matrix 
$\bfTheta \in \reals^{N_B \times N}$,
where $N_B$ is the number of grid points where the deformation $\DB$ is
computed (in our experiments we use $N_B = 1001$ over an interval that
extends 100 km on either side of the fault region).  The Okada model
implemented in the GeoClaw {\tt dtopotools} Python module is used, which
directly computes $\DB$ from $\bs$ and so we do not actually compute the
matrix $\bfTheta$, but it is useful conceptually.  In particular, if the K-L
expansion \cref{sm} is to be used to compute $\bs$ then we find that
\begin{equation}\label{DBm}
\DB = \bfTheta \bmu + \bfTheta \bfV \bfLambda^{1/2} \zm.
\end{equation} 
The vector $\bfTheta\bmu$ is obtained by applying Okada to the mean slip
vector.  The matrix $\bfTheta \bfV$ can be computed by applying Okada to each
column of $\bfV$  to compute the columns of the product matrix.
Since the sum only involves $m$ nonzero terms, we need only apply Okada to
columns 1 through $m$ of $\bfV$ (i.e. the first $m$ K-L modes $\bv_1,~\bv_2,~
\ldots,~\bv_m$  used to express $s$).  Hence if we plan
to use at most 20 modes of the K-L expansion then we need only apply the
Okada model to 21 slip vectors and we can then take linear combinations of the
resulting sea floor deformations, rather than applying Okada to 20,000 slip
realizations separately. 

In practice this can be simplified further.  Applying Okada to a mode $\bv_k$
actually requires applying Okada to each of the $N$ subfaults, weighting by the
corresponding element of $\bv_k$, and summing over all the subfaults.  So
applying Okada to $m$ modes in this way actually requires applying Okada $mN$
times.   Instead, 
we can first apply the Okada model to $N$ unit source scenarios in which the
slip is set to 1 on the $j$th subfault and to 0 on all other subfaults.  Call
this slip vector $\bs^{[j]}$.  Applying Okada to this gives a resulting
$\DB^{[j]} = \bfTheta \bs^{[j]}$.  Now for any slip vector $\bs$ we can compute
$\bfTheta\bs$ as
\begin{equation}\label{DBunit}
\bfTheta\bs = \sum_{j=1}^N s_j \DB^{[j]}.
\end{equation} 
In particular taking $\bs = \bv_k$ would give $\bfTheta \bv_k$, but \cref{DBunit} can
be used directly to compute the seafloor deformation $\DB = \bfTheta\bs$ for any
slip realization.  This approach can also be used in the lognormal case
described in \Sec{sec:lognormal} and employed in \Sec{sec:def2d}.

{\bf Subsidence or uplift.}
One quantity that has a significant impact on the severity of tsunami
flooding is the vertical displacement of the seafloor at the coast.  
If this displacement
is negative and the land subsides at the shore, then flooding will generally
be much worse than if uplift occurs.  The behavior seen for a particular
event depends on how far offshore the subduction zone is, which is generally
directly related to the width of the continental shelf offshore from the
community of interest (since the top edge of the fault is usually located near
the trench at the edge of the shelf, which is part of the continental 
plate beneath which the oceanic plate is subducting). This distance can vary
from only a few km (e.g. along the Mexico coast) to 200 km (e.g. along the
coast of Japan where the Tohoku event occurred).    In our model the top
of the plate is at $x=0$ and we choose the coast line location to be at
$x=75$ km, which gives a wide range of subsidence and uplift values,
as can be observed for the realizations shown in \Fig{fig:realiz1}.

\begin{figure}
\hfil\includegraphics[width=0.45\textwidth]{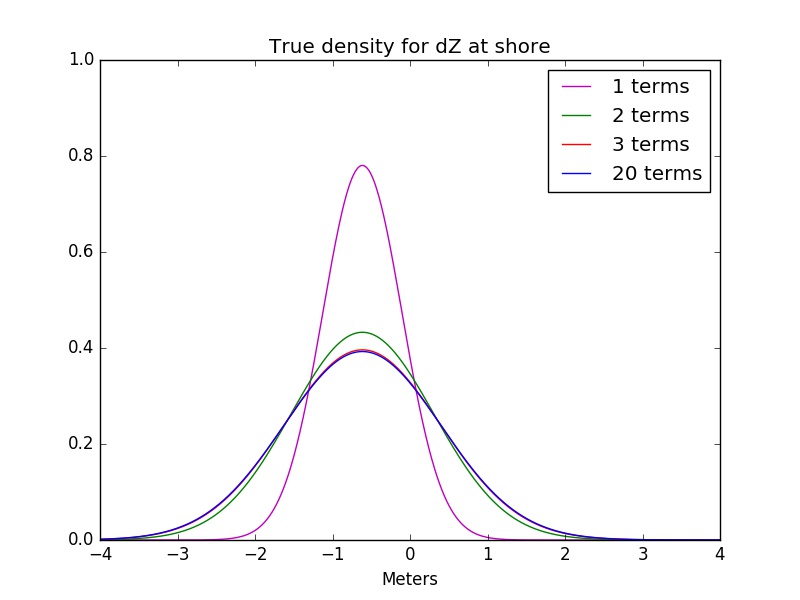}\hfil
\hfil\includegraphics[width=0.45\textwidth]{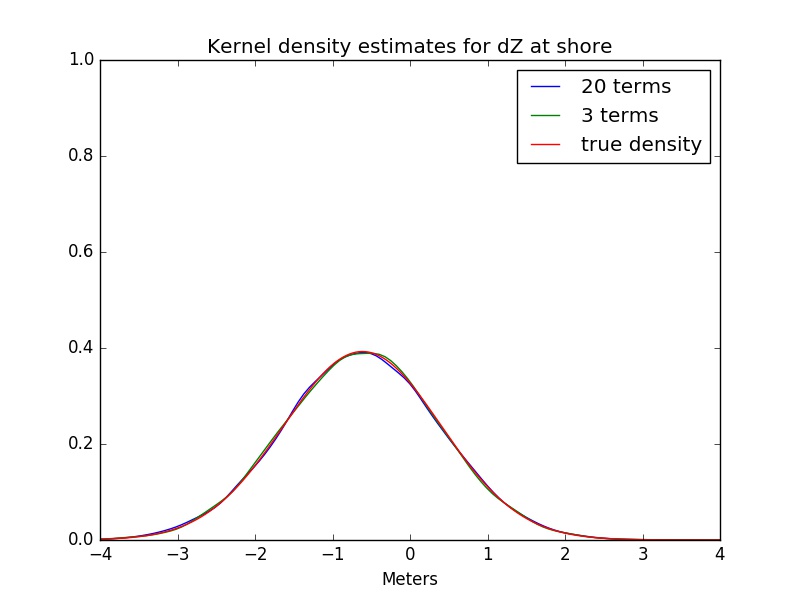}\hfil
\caption{\label{fig:DBshore-density} 
The figure on the left shows the true Gaussian density for the shoreline
displacement for K-L expansions with 1, 2, 3, or 20 terms. 
The figure on the right shows the kernel density estimate from 20,000
samples using 3 terms or 20 terms, together with the true density for 20
terms.  
  }
\end{figure}

The displacement at one particular point is easy to determine from each
realization, it is just one entry in the vector of sea floor deformation
obtained from the Okada model, say $\DB_j = \be_j^T \DB$ 
for some $j$, where $\be_j$ is the unit vector with a 1 in position $j$.
(We have evaluated
$\DB$ on a fine grid so we assume we don't need to interpolate).  
As such, this particular quantity is in fact easy to compute directly from
$\bz$ for any given realization, as
\begin{equation}\label{DBshore}
\DBshore = \DBshore(\mu) + \bb^T\bz,
\end{equation} 
where $\DBshore(\mu) = \be_j^T \bfTheta\bmu$ 
is the shoreline displacement resulting from the mean
slip and the row vector $\bb^T$ is
\begin{equation} \label{DBb}
\bb^T = \be_j^T\bfTheta \bfV \bfLambda^{1/2},
\end{equation} 
i.e. the vector consisting of the $j$th component of the sea floor
displacement resulting from applying Okada to each K-L mode, scaled by the
square root of the corresponding eigenvalue.
From \cref{DBshore} it follows immediately that $\DBshore$ is normally
distributed with mean $\DBshore(\mu)$ and variance $\sigma^2 = \sum_{k=1}^m
b_k^2$ (in the Gaussian case considered here, not in the lognormal case
considered below).  Hence for this particular quantity of interest in the 
Gaussian case, we do not need to
estimate the statistics based on a large number of samples. We can
immediately plot the Gaussian density function for the ``full'' expansion with 20
terms and compare it to the density for the truncated series with only 3
terms.  These are seen to lie nearly on top of one another in
\Fig{fig:DBshore-density}.  This plot also shows the density that would be
obtained with only 1 or 2 terms in the K-L expansion, which are
substantially different.  
This confirms that, in terms of this particular quantity
of interest, it is sufficient to use a 3-term K-L expansion (but not fewer
terms).

\Fig{fig:DBshore-density} also shows the density as estimated using 20,000
samples, using a kernel density estimate
computed using the Python package {\tt seaborn} of \cite{seaborn}.  
With either 3 terms or 20 terms, the estimated density lies nearly on top of
the true density, giving confidence  
that the sampling has been programmed properly and that 20,000
samples is sufficient since the true density is known in this case.

\begin{figure}
\hfil\includegraphics[width=.45\textwidth]{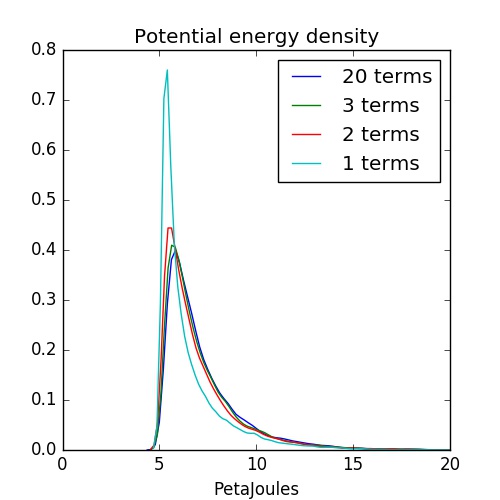}\hfil
\hfil\includegraphics[width=.45\textwidth]{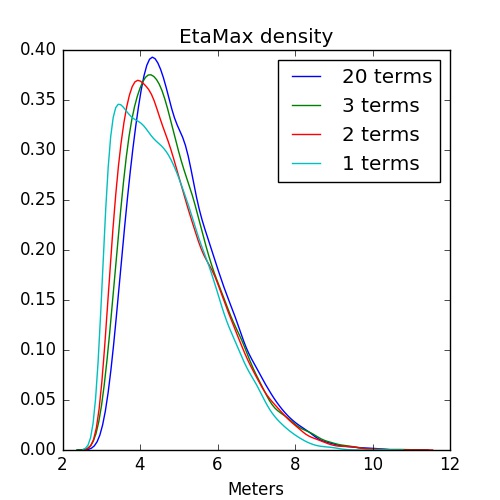}\hfil
\caption{\label{fig:kde} 
Kernel density estimates based on 20,000 samples, using 1, 2, 3, or 20 terms
in the K-L expansion.  The left figure shows the potential energy \cref{PE}
and the right figure shows the maximum amplitude of deformation (sea surface
elevation).
  }
\end{figure}

{\bf Potential energy.}
From the samples it is possible to also estimate the densities for other
quantities of interest for which it is not possible to compute the true
density.  We consider two additional quantities that have relevance to the
magnitude of the tsunami generated.  One is the potential energy of the
initial perturbation of the ocean surface, which is one measure of its
potential for destruction.  The potential energy is given by
\begin{equation}\label{PE}
E = \frac 1 2 \int\int \rho g \eta^2(x,y)\, dx\,dy
\end{equation} 
where $\rho = 1000$ kg/m$^3$ is the density of water, 
$g = 9.81$ m/s$^2$, and $\eta(x,y)$ is
the initial perturbation of the surface.  With our assumption that the sea
surface moves instantaneously with sea floor deformation generated from the
slip, $\eta$ is equal to the sea floor displacement in the ocean, 
while onshore we set  $\eta=0$
since displacement at these points does not contribute to the potential
energy of the tsunami.  For the one-dimensional problem considered here, we
sum the square of the displacement over $x < 75$ km and scale
by $\rho g L \Delta x$ to define $E$, taking $L=100$ km.  
Finally we multiply by $10^{-15}$ so that the results are order 1, with
units of PetaJoules.
We use the kernel density estimator of
{\tt seaborn} to plot the results obtained with 20,000 samples, using 20
terms or truncating further to 1, 2, or 3 terms.  The results in 
\Fig{fig:kde} again show that 3 terms is sufficient to obtain
very similar results to 20 terms.

{\bf Maximum wave height.}
The maximum positive seafloor displacement gives the maximum amplitude of the
tsunami at the initial time.  We expect this to be
positively correlated with the amplitude of the wave that approaches shore
(although the wave propagation can be complicated by multiple peaks, the
location of the $\etamax$ relative to the shore, or various other factors
that can only be studied with a full tsunami model).
The right plot of \Fig{fig:kde} shows the kernel density estimates
of this quantity $\etamax$.

{\bf Joint probability densities.}
It is also interesting to plot the joint
probability density of pairs of quantities to better explore the ability of 3
terms to capture the variation.
This is illustrated in \Fig{fig:joint}, where the top row shows kernel
density estimates (again computed using {\tt seaborn}) for $E$ vs.\ $\etamax$
and the bottom rows shows $\DBshore$ vs.\ $\etamax$.  In each case the left figure
shows the density computed from 20,000 realizations of the 20-term K-L
expansion while the right figure shows the density estimated from an equal
number of 3-term expansions.  In each case it appears that the
3-term expansion captures the bulk of the variation.

The joint distribution of $\etamax$ and $\DBshore$
is of particular interest since the most dangerous
events might be those for which $\etamax$ is large while $\DBshore$ is most
negative (greatest subsidence of the coast).  The fact that the joint
distributions look quite similar gives hope that the 3-term model will
adequately capture this possibility.

\begin{figure}
\hfil\includegraphics[width=.9\textwidth]{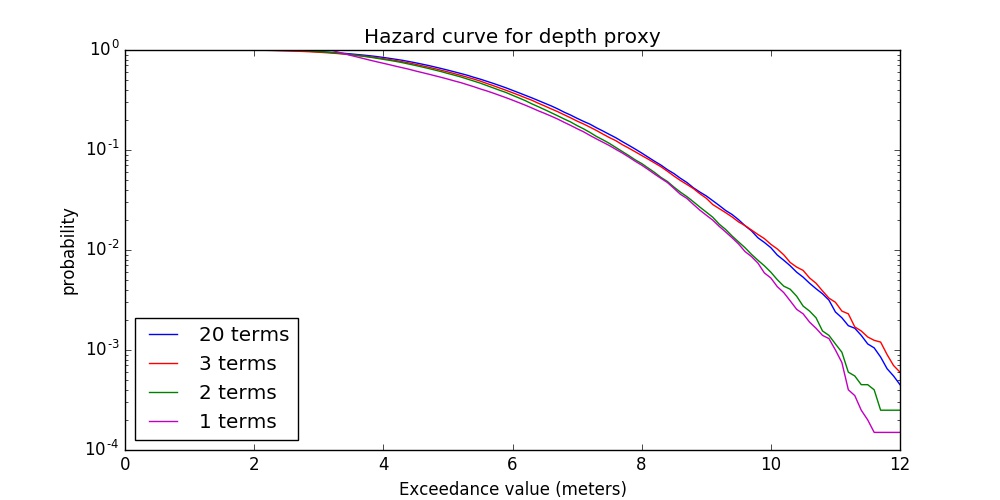}\hfil
\caption{\label{fig:hazcurves1} 
Hazard curves based on the proxy for flooding depth given by $\etamax -
\DBshore$.  Based on 20,000 samples using the full
20-term K-L expansion, compared with the hazard curves obtained using only 1,
2, or 3 terms in the expansion. Note that 3 terms is sufficient to obtain the
hazard curve to high precision.
  }
\end{figure}

{\bf Depth proxy hazard curves.}
The goal of a full-scale PTHA exercise is often to generate hazard curves at
many points onshore or in a harbor.  A hazard curve shows, for example, the
probability that the maximum flow depth will exceed some value as a function
of that exceedance value.  Construction of these curves is discussed, for
example, in the appendices of \cite{GonzalezLeVequeEtAl2014}.
The curves may vary greatly with spatial location due
to the elevation of the point relative to sea level, 
and also due to the manner in which a tsunami
interacts with nearby topography.  Hazard curves must thus be computed using
fine-grid simulations of the tsunami dynamics and cannot be computed
directly from the sea floor deformation alone.  However, as a proxy for
flooding depth we might use $D = \etamax - \DBshore$, the maximum offshore sea
surface elevation augmented by any subsidence that occurs at the shore.  We
do not claim that this is a good estimate of the actual maximum water depth
that will be observed at the shore, but computing hazard curves for this
quantity provides another test of how well the 3-term K-L expansion captures the
full probability distribution described by the 20-term expansion.
This curve is obtained by computing
$D$ for each sample and determining the fraction of samples
for which this is above $\zeta_i$, for each exceedance level $\zeta_i$ on a
fine grid covering the range of $D$ observed.
\Fig{fig:hazcurves1} shows the resulting hazard curve obtained with the 20-term
expansion.  The curve obtained with the 3-term expansion is also shown, and
lies nearly on top of it.

\begin{figure}
\hfil\includegraphics[width=0.45\textwidth]{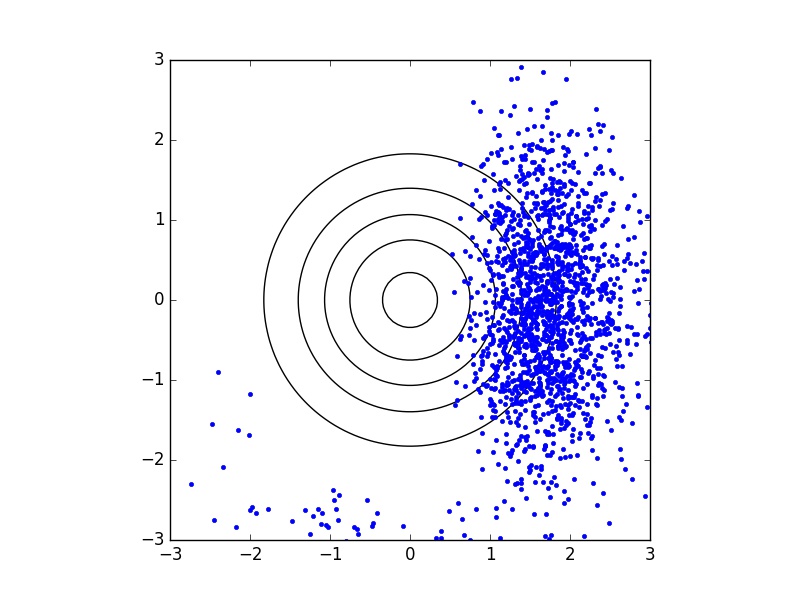}\hfil
\hfil\includegraphics[width=0.45\textwidth]{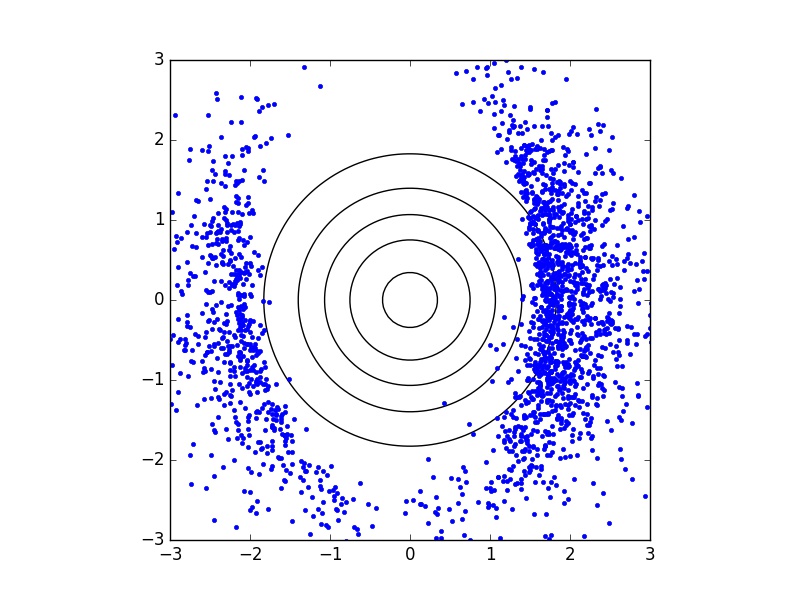}\hfil
\caption{\label{fig:bigz1z2} 
Scatter plots in the $z_1$-$z_2$ plane of the subset of 3-term events for
which the proxy depth is greater than 8 m (left) or  for which 
the potential energy is greater than 9.5 PetaJoules (right).
  }
\end{figure}

{\bf Exploring parameter space.}
One advantage of describing the probability space of possible events in
terms of a small number of stochastic parameters is that it may be possible
to examine structures in this space, which can be important in developing a
cheap surrogate model to use in estimating probabilities and computing
hazard curves for practical quantities of interest.  
For example, we can ask what parts of parameter space lead to the worst
events. The left figure in
\Fig{fig:bigz1z2} shows the events (projected to the $z_1$-$z_2$ plane)
from the above tests with the 3-term K-L expansion
for which the proxy depth is greater than 8 m.  The
contours of the bivariate normal distribution are also plotted. A scatter
plot of all 20,000 events would cluster in the middle, but we observe that the
events giving this extreme depth tend to have $z_1 > 1$.  From \Fig{fig:modes1d}
we see that positive $z_1$ redistributes slip from the down-dip to the up-dip 
portion of the fault.  This agrees with common wisdom from past events that
concentration near the up-dip edge gives particularly severe tsunamis (as in
the case of the 2011 Tohoku event).  
The right figure in \Fig{fig:bigz1z2} shows a similar scatter plot 
of $z_1$-$z_2$ values for which the potential energy was above 9.5 PetaJoules.
In this case most of the extreme events have $z_1$ either very
positive or very negative.  In the latter case slip is concentrated toward 
the down-dip portion of the fault, which leads to a smaller maximum surface
displacement but the displacement spreads out further spatially for a 
deep rupture which can lead to large potential energy since this is
integrated over space.

\section{Lognormally distributed slip}\label{sec:lognormal}
If we wish to instead generate slip realizations 
that have a joint lognormal distribution
with a desired mean and covariance matrix, we can first generate realizations of a
joint Gaussian random field and then exponentiate each component.  
This approach will be used in the two-dimensional example below in
\Sec{sec:def2d}.

In this case we first choose the desired mean $\bmu$ and covariance matrix
$\bsigma$ for the slip,  and then compute the necessary mean $\bmu^g$ and
covariance matrix $\bsigma^g$ for the Gaussian to be generated by the K-L
expansion, using the fact that if $\bg$ is a random variable from ${\cal
N}(\bmu^g, \bsigma^g)$, then $\exp(\bg)$ is lognormal with mean and covariance
matrix given by:
\begin{equation}\label{logn}
\mu_i = \exp(\mu_i^g + \hat C_{ii}^g/2), \qquad \hat C_{ij} =
\mu_i\mu_j(\exp(\hat C_{ij}^g) - 1).
\end{equation} 
Hence we can solve for
\begin{equation}\label{logng}
\begin{split} 
\hat C_{ij}^g &= \log(\hat C_{ij}/\mu_i\mu_j + 1),\\
\mu_i^g &= \log(\mu_i) - \frac 1 2 \hat C_{ii}^g.
\end{split}
\end{equation} 
We now find the eigenvalues $\lambda_k$ and eigenvectors $\bv_k$
of $\bsigma^g$.  To generate a realization we choose $N$ values $z_k\sim{\cal
N}(0,1)$ and then form the K-L sum
\begin{equation}\label{KL1g}
\bs^g = \bmu^g + \sum_{k=0}^N z_k \sqrt{\lambda_k} \bv_k.
\end{equation} 
We then exponentiate each component of $\bs^g$ to obtain the 
slip values, which then have the desired joint lognormal distribution
(see e.g., \cite{Ghanem1999}).

As described, this will generate realizations with total slip (and hence
magnitude $M_w$) that vary around the mean.  As in the Gaussian case, we can
drop the nearly-constant $\bv_0$ term from the sum to reduce this variation.  
We can also generally truncate the series to a much smaller number of terms 
and still capture most of the variation if the eigenvalues are rapidly
decaying.

Now consider the special case where 
we make the same assumptions as in \Sec{sec:KL} that
$\hat C_{ij} = \sigma_i\sigma_jC_{ij}$ where $\bfC$ is the desired 
correlation matrix and
$\sigma_i = \alpha \mu_i$, while the mean $\mu_i$ was given by some taper
$\tau_i$ scaled by a scalar value $\bar \mu$.  
Then computing
$\bmu^g$ and $\bsigma^g$ according to \cref{logng}, we find that:
\begin{equation}\label{logng2}
\begin{split} 
\hat C_{ij}^g &= \log(\alpha^2 C_{ij} + 1),\\
\mu_i^g &= \log(\bar\mu\tau_i) - \frac 1 2 \log(\alpha^2 + 1)
\end{split}
\end{equation}
We see that the covariance matrix in this case depends only on the correlation
matrix and the scalar $\alpha$, not on the mean slip itself (and in particular
is independent of the taper).  We also find that
$\exp(\mu_i^g) = \bar\mu\tau_i / \sqrt{\alpha^2+1}$ is simply a scalar multiple
of the taper.

Using these assumptions and the fact that
\[
\exp\left(\bmu^g + \sum_{k=1}^N z_k \sqrt{\lambda_k} \bv_k\right)
= \exp(\bmu^g) \exp\left(\sum_{k=1}^N z_k \sqrt{\lambda_k} \bv_k\right),
\]
it is easy to generate realizations that have exactly the desired
magnitude: simply compute 
\begin{equation}\label{expsum}
\exp\left(\sum_{k=1}^N z_k \sqrt{\lambda_k} \bv_k\right),
\end{equation} 
multiply the result by the desired taper, and then rescale by a multiplicative
factor so that the area-weighted sum of the slips gives the total slip
required for the desired seismic moment.

\section{Two-dimensional case}\label{sec:def2d}

\begin{figure}[t]
\hfil\includegraphics[width=.9\textwidth]{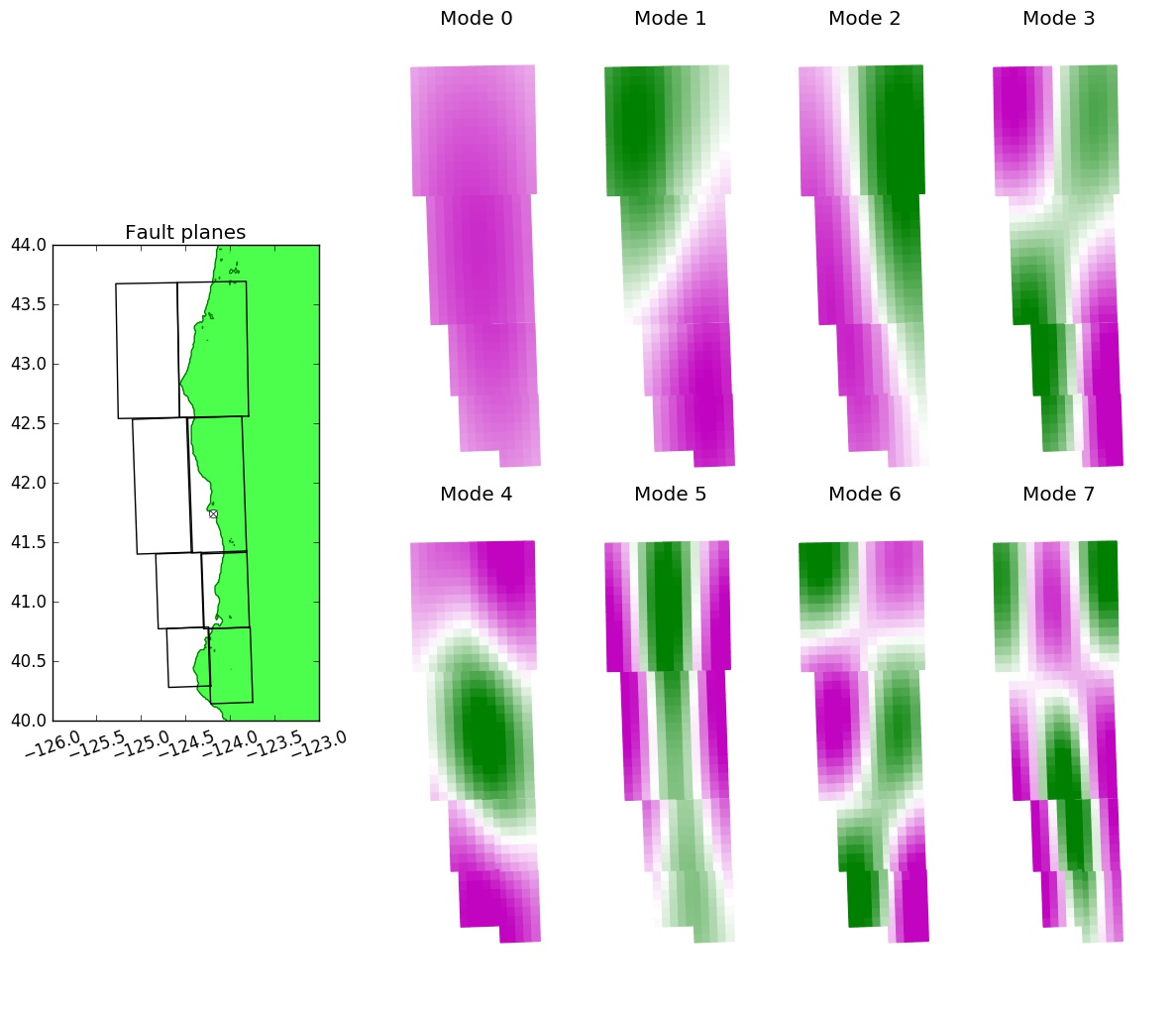}\hfil

\caption{\label{fig:csz2} 
Southern portion of the CSZ fault showing location of Crescent City, CA and the
8 subfaults that are further subdivided into 540 subfaults.  The first 7
eigenmodes of the resulting covariance matrix $\bsigma^g$ are also shown.
  }
\end{figure}

We now present an example in which the slip is allowed to vary
in both directions along a fault surface.  For illustration we use a subset of
the Cascadia Subduction Zone from \Fig{fig:csz1}, taking only the
southern-most 8 fault segments, as illustrated in \Fig{fig:csz2}. 
These are subdivided into 540 smaller fault planes for the purposes of defining
the slip.  

To define the $540 \times 540$ correlation matrix, we need to compute the
pairwise ``distance'' between subfault $i$ and subfault $j$.  
We can compute the Euclidean distance $d_{ij}$, but
for this fault geometry we expect
a longer correlation length in the strike direction than down-dip, 
so we wish to define
\begin{equation}\label{corr2d}
C_{ij} = \exp(-(\dstrike(i,j)/\rstrike) - (\ddip(i,j)/\rdip))
\end{equation} 
where $\dstrike(i,j)$ and $\ddip(i,j)$ are estimates of the distance between
subfaults $i$ and $j$ in the strike and dip direction, respectively, and
$\rstrike,\rdip$ are the correlation lengths in each direction.  We define
$\ddip(i,j)$ using the difference in depth between the two subfaults 
and the dip angle $\delta$
as $\ddip(i,j) = \ddepth / \sin(\delta)$,
setting $\dstrike(i,j) = \sqrt{d_{ij}^2 - \ddip(i,j)^2}$.
We take the correlation lengths to be 40\% of the fault length and width
respectively, $\rstrike = 130$km and $\rdip = 40$km.
We again use an exponential autocorrelation function as defined in
\cref{corr2d}, but this could easily be replaced by a different ACF.  
We use the lognormal approach described in \Sec{sec:lognormal}, with parameter 
$\alpha = 0.5$.  
\Fig{fig:csz2} shows the first 8 eigenmodes of $\bsigma^g$.
Again we drop Mode 0 from the sum, since this mode is roughly constant over the
fault.

\begin{figure}[t]
\hfil\includegraphics[width=.9\textwidth]{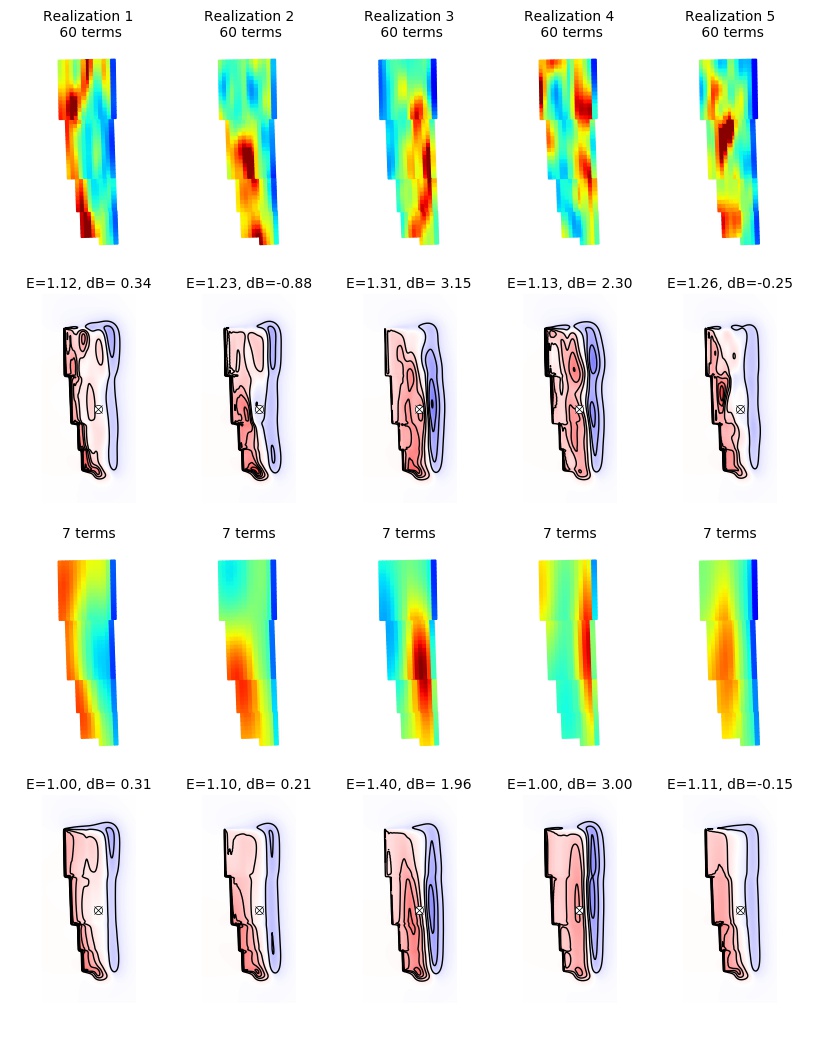}\hfil

\caption{\label{fig:CSZrealizations} 
The top row shows 5 sample realizations of slip on the southern CSZ fault,
as computed with a 60-term K-L expansion. 
The second row shows the resulting seafloor deformation, with an indication of
the potential energy and the vertical displacement at Crescent City, CA, which
is indicated by the X in the figures.  The third row shows the same 5
realizations but with the K-L series truncated to 7 terms, and the bottom row
shows the resulting seafloor deformations.
  }
\end{figure}

To create slip realizations, we use \cref{expsum} and then apply 
a tapering only at the
down-dip edge, given by \cref{taper} with $d_{max}=20000$m.
We then scale the slip so that the resulting seismic moment gives $M_w=8.8$.
\Fig{fig:CSZrealizations} shows 5 typical realizations,
comparing the slip generated by a 60-term K-L expansion with the slip generated
when the series is truncated after 7 terms.  The resulting seafloor
deformation in each case is also shown, along with the potential energy and the
subsidence/uplift $\DBshore$
at one point on the coast, the location of Crescent City, CA.
Note that in each case the 7-term series gives a smoother version of the slip
obtained with 60 terms, and the seafloor deformations are more similar than the
slip patterns, as expected from the one-dimensional analogous case shown in
\Fig{fig:realiz1}.  The potential energy and $\DBshore$ are also seen to be
similar when the truncated series is used to the values obtained with the
longer 60-term series.

\begin{figure}
\hfil\includegraphics[width=.45\textwidth]{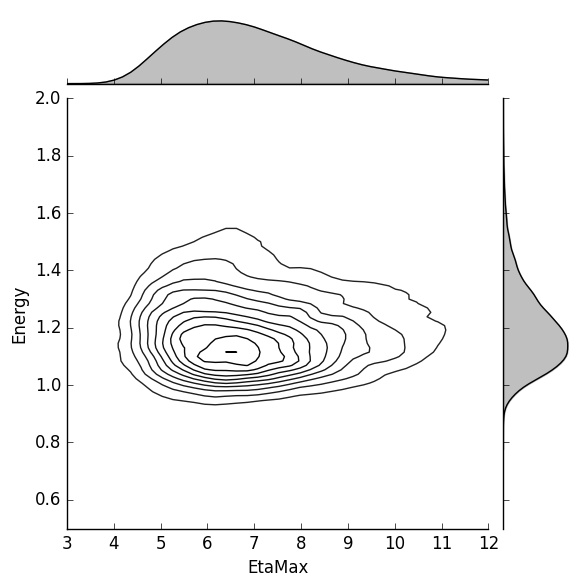}\hfil
\hfil\includegraphics[width=.45\textwidth]{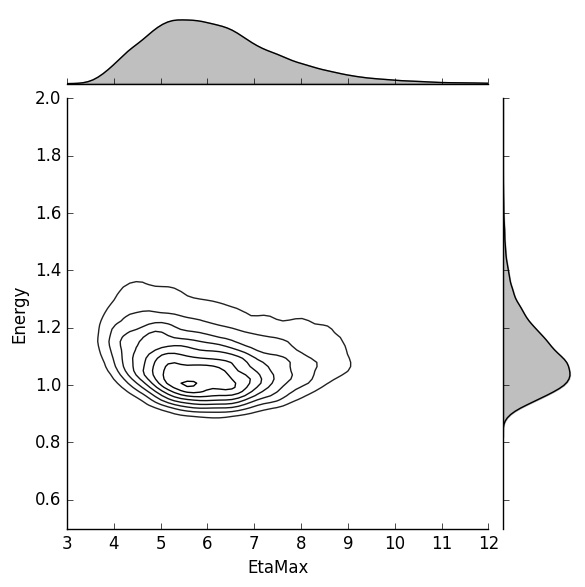}\hfil

\hfil\includegraphics[width=.45\textwidth]{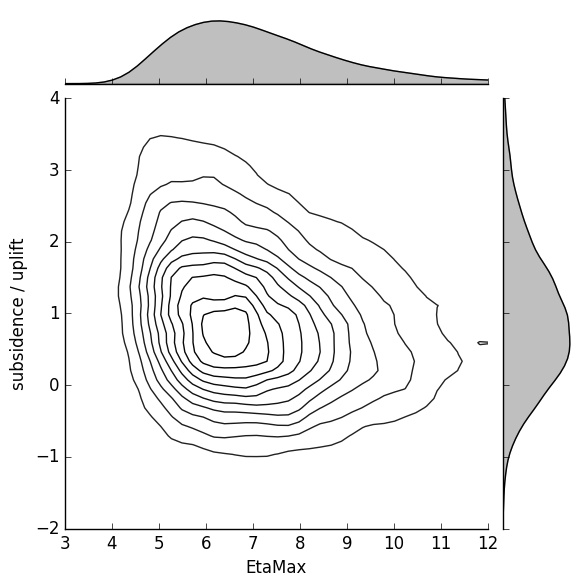}\hfil
\hfil\includegraphics[width=.45\textwidth]{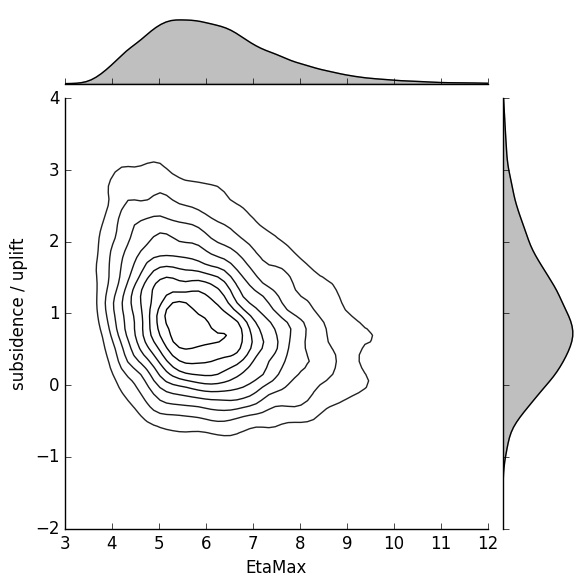}\hfil

\caption{\label{fig:joint2} 
Joint and marginal probability densities for different quantities, comparing
the densities estimated using the 60-term expansion (left column) and the
7-term expansion (right column) for the two-dimensional fault case.  
The top row shows the joint density of $\etamax$ with
potential energy $E$ of the tsunami generated.
The bottom row shows the joint density
of $\etamax$ with $\DBshore$, the vertical displacement at Crescent City, CA.
In each case 20,000 realizations similar to those shown in
\Fig{fig:CSZrealizations} were used to create these kernel density estimates.
  }
\end{figure}

We can explore the statistical properties by repeating any of the experiments
performed above in the one-dimensional case.  In the interest of space, we only
show one set of results, the same joint and marginal densities examined in the
one-dimensional case in \Fig{fig:joint}.  The comparisons for the
two-dimensional fault are shown in \Fig{fig:joint2}.
To generate each column of figures we computed 20,000 slip realizations and the
resulting seafloor deformations (via \cref{DBunit}). The first column shows
statistics when a 60-term KL-expansions is used, producing realizations
similar to those shown in the top row of \Fig{fig:CSZrealizations}.
The second column of figures was produced using an independent set of 7-term
realizations (i.e. these were not obtained by truncating the 60-term series
from the first set, but rather by generating 20,000 independent samples).
Even in this two-dimensional case, less than
10 minutes of CPU time on a MacBook Pro laptop was required to generate each
set of 20,000 realizations, the resulting seafloor deformations,
and the kernel density plots.

\section{Discussion}\label{sec:discussion}
We have presented an approach to defining a probability distribution for
earthquake slip patterns on a specified fault geometry that has been subdivided
into an arbitrary number of rectangular subfaults, with a specified mean
and covariance matrix.  Slip realizations can be generated that 
either have a joint normal distribution
or a joint lognormal distribution.  Once the parameters have been chosen that
define the distribution, it is very easy to generate an arbitrary number of
sample realizations from the distribution, simply by drawing the coefficients
$z_k$ of the K-L series from independent normal distributions.

We have also illustrated that with a realistic choice of correlation length,
the K-L series can be truncated to a relatively small number of terms.  For
tsunami modeling applications, the Okada model is applied to each slip pattern
to generate the seafloor deformation and it was shown that this is a smoothing
operation that can further reduce the number of terms needed, and hence the
dimension of the stochastic space that must be explored in doing PTHA analysis.

To use this approach for practical PTHA analysis, two major challenges must be
addressed.  The first is to tackle the {\em epistemic uncertainty} associated
with the lack of knowledge about possible future earthquakes.  We would like to
choose the parameters defining the probability
distribution in a suitable way for real fault geometries so that it accurately
represents the space of possible future earthquakes.  
However, realistic specification of these critical seismic parameters and
quantifying the associated uncertainties and geophysical contraints is a
major challenge that \cite{SteinGellerEtAl2012}
have reviewed and summarized; they
characterize the problem as a failure of earthquake hazard mapping, in
general, and make recommendations regarding improvements.  The problem is
particularly severe in the case of near-field PTHA studies, because tsunami
impact on a coastal community is highly sensitive to details of the seismic
deformation (e.g., \cite{Geist2002}).
Existing expertise and geophysical constraints should at least
be incorporated in the
choice of these parameters.  The ability to generate many realizations and
examine statistics of quantities such as those used in this paper may help in
this.  As one example, the parameters chosen in this paper for the CSZ example
tend to give uplift rather than subsidence at Crescent City (as can be seen in
the marginal distribution of subsidence/uplift in \Fig{fig:joint2}).  If this
is viewed as inconsistent with the geological evidence from past events, this
could be adjusted, for example by tapering the slip more on the down-dip side.
Moving more of the slip up-dip will cause more subsidence at the shore.  
It would also be possible to explore ways in which the epistemic
uncertainty associated with the lack of knowledge about the true probability
distribution affect the resulting hazard maps generated by a PTHA analysis, for
example by doing the analysis with different parameter choices, and hence
different probability distributions, to see how robust the PTHA analysis is to
changes in assumptions.

The second major challenge is to deal with the {\em aleatoric uncertainty} that
is still present even if the parameters defining the probability distribution
were known to be correct.  We are still faced with a high-dimensional space to
sample in order to perform PTHA analysis.  For example, if we wish to compute a
hazard curve similar to \Fig{fig:hazcurves1} for the probability that the
maximum depth $D$ at some particular point
will exceed various depths, then for each exceedence value
$D_e$  we need to calculate 
\begin{equation}\label{excprob}
P[D > D_e] = \int \rho(\bz)I(\bz;D_e)\,d\bz
\end{equation} 
where the integral is over the $m$-dimensional stochastic space of coefficients
$z$ of the K-L sum (assuming $m$ terms are used) and
$I(\bz;D_e)$ is an indicator function that is
1 at points $\bz\in\reals^m$ where the corresponding realization gives
a tsunami that exceeds $D_e$ and 0 elsewhere (or it could take values between
0 and 1 to incorporate other uncertainties, e.g. if the approach of
\cite{Adamstide} is used to incorporate tidal uncertainty). The function
$\rho(\bz)$ in \cref{excprob} is the probability density for $\bz$.  
In the K-L approach, 
$\bz$ is a vector of i.i.d. normally distributed values so $\rho(\bz)$ 
is known; for the $m$-term expansion it takes the form
\begin{equation}\label{rhogauss}
\rho(\bz) = \frac{1}{\sqrt{(2\pi)^m}} \exp\left(-\frac 1 2 \sum_{i=1}^m
z_i^2\right).
\end{equation} 
To generate
\Fig{fig:hazcurves1}, we used a simple Monte-Carlo method in which the
integral in \cref{excprob} is replaced by
\begin{equation}\label{mc}
P[D > D_e]\approx \frac 1 {n_s} \sum_{j=1}^{n_s} I_e(\bz^{[j]}),
\end{equation} 
with $n_s=$ 20,000 samples and $\bz^{[j]}$ now represents the 
$j$th sample, drawn
from the joint normal distribution with density $\rho(\bz)$.  
This is feasible with the depth proxy used here,
but would not be possible if a full tsunami model is used to compute $D$, which
may take hours of computing time for each sample.  
The number of simulations required can be reduced by source-filtering
techniques that identify a ``most-important'' subset of  realizations that
contribute most to the tsunami impact on a particular site, e.g.
\cite{LoritoSelvaEtAl2015}.
An alternative would be to
compute the integral with a quadrature algorithm based on sampling on a grid in
$z$-space, but this is infeasible for high dimensions $m$. For example, if
$m=10$ then a tensor-product 
grid with only 4 points in each direction has $4^{10}\approx 10^6$ points.

Many other techniques have been developed in recent years to estimate such
integrals in high dimensional spaces, including for example
Latin hypercube sampling (e.g., \cite{OlssonSandberg2002})),
sparse grids (e.g., \cite{NobileTemponeWebster2008}), and
quasi-random grids (e.g., \cite{DickKuoSloan2013}) 
that have many fewer points than uniform
tensor-product grids.  There are also several Monte-Carlo
sampling methods that can obtain accurate results with many fewer samples than
the naive sum of \cref{mc}, including multi-level or multi-fidelity methods
(e.g., \cite{Cliffe2011,Giles2008,Peherstorferetal2016})
that combine results from many simulations that are cheap to compute 
with a relatively few simulations with the full model on a fine grid.  
Cheaper approximations might be obtained by
using some of the proxy quantities from this paper, by computing with a full
tsunami model but on coarse grids, or by developing surrogate models or
statistical emulators based on relatively few samples (e.g.
\cite{BastosOHagan2009,benner15,LiLiXiu2011,SarriGuillasEtAl2012}).
We are currently exploring several of these approaches for PTHA and will report
on them in future publications.

The computer code used to generate the examples used in this paper 
can be found in the repository \url{https://github.com/rjleveque/KLslip-paper}
and will be archived with a permanent DOI once this paper is finalized.

\vskip 10pt
{\bf Acknowledgements.}
The initial phase of this work was performed when GL was employed at Pacific
Northwest National Laboratory (PNNL) and KW was a postdoctoral fellow
supported in part by PNNL and by the University of Washington.
This work was also supported in part by NSF grants DMS-1216732 and EAR-1331412, 
and funding from FEMA. 
The authors have benefitted from discussions with many applied mathematicians
and geoscientists concerning the approach developed in this paper,
including in particular Art Frankel, Finn L{\o}vholt, Martin Mai, Siddhartha
Mishra, Diego Melgar, and Hong Kie Thio.

\vskip 1cm
\bibliographystyle{spbasic}

\end{document}